\documentclass[a4paper,10pt]{article}
\usepackage[utf8x]{inputenc}
 \usepackage[T1]{fontenc} 
 \usepackage[finnish,english]{babel}

\usepackage{amsmath}
\usepackage{amsthm}
\usepackage{amssymb}
\usepackage{mathtools}

\newcommand{\integers}{\mathbb{N}}

\title{Surface subgroups on 1-vertex and 3-vertex polyhedra forming triangular hyperbolic buildings}
\author{Riikka Kangaslampi}

\begin{document}

\maketitle

\begin{abstract}
In this article we study surface subgroups of groups acting simply transitively on vertex sets of certain hyperbolic triangular buildings. We concentrate on surface subgroups on 1-vertex and 3-vertex polyhedra that can be used to define these buildings. In \cite{KV2} and \cite{CKV} the authors constructed and classified all groups acting simply transitively on the vertices of hyperbolic triangular buildings of the smallest non-trivial thickness. These groups gave the first examples of cocompact lattices acting simply transitively on vertices of hyperbolic triangular Kac-Moody buildings that are not right-angled. They constructed the groups with the polygonal presentation method introduced in \cite{Vdovina}. The hyperbolic buildings arise as universal covers of polyhedra glued together from with 15 triangular faces with words written on the boundary.  In \cite{KV3} the authors proved, that most of the obtained 23 torsion free groups do not admit periodic planes of genus 2. Here we take another approach to study surface subgroups of these groups. Namely, we consider first 2-cycles in the 1-vertex polyhedron defined by the triangles, then in the 3-vertex cover of this polyhedron. As a result we find surface subgroups only in three of the 23 torsion free groups. These were also found in \cite{KV3}.

\end{abstract}

\section{Introduction}

In \cite{KV2} the authors classified all torsion-free groups acting simply transitively on the vertices of hyperbolic triangular buildings of the smallest non-trivial thickness. They constructed the groups with the polygonal presentation method introduced in \cite{Vdovina}. As a result, they obtain 23 non-isomorphic groups, each defined by 15  generators $x_1,x_2,\ldots, x_{15}$ and 15 cyclic relations, each of them of the form $x_ix_jx_k=1$, where not all the indices are the same. The underlying hyperbolic building is the universal cover of the polyhedron glued together from 15 geodesic triangles with angles $\pi/4$ and with the letters from the relations written on the boundary. In constructing the polyhedron the sides of the triangles with the same labels are glued together, respecting the orientation. For example, the presentation $T_1$ 
obtained in \cite{KV2} are given in Table \ref{T1T21}. 

\begin{table}[h!]\label{T1T21}
\centering
    \begin{tabular}{l}
$(x_1,x_1,x_{10})$\\
$(x_1,x_{15},x_2)$\\
$(x_2,x_{11},x_9)$\\
$(x_2,x_{14},x_3)$\\
$(x_3,x_7,x_4)$\\
$(x_3,x_{15},x_{13})$\\
$(x_4,x_8,x_6)$\\
$(x_4,x_{12},x_{11})$\\
$(x_5,x_5,x_8)$\\ 
$(x_5,x_{10},x_{12})$\\
$(x_6,x_6,x_{14})$\\
$(x_7,x_7,x_{12})$\\
$(x_8,x_{13},x_{9})$\\ 
$(x_9,x_{14},x_{15})$\\
$(x_{10},x_{13},x_{11})$\\

    \end{tabular}
      \caption{Triangles in presentation $T_1$  from \cite{KV2}.}
    \end{table}

Thus these sets of 15 triangles, with angles $\pi/4$, words specified in \cite{KV2} written at the boundary and glued together respecting orientation, all give a polyhedron that has one vertex and the smallest generalised quadrangle as the link.  The universal cover of this polyhedron is a hyperbolic triangular building, and the group with 15 generators $x_1,x_2,\ldots, x_{15}$ and the 15 words from the boundaries of the triangles as relations, acts on the building cocompactly and simply transitively. These groups are the first examples of cocompact lattices acting simply transitively on vertices of hyperbolic triangular Kac-Moody buildings that are not right-angled. For a general introduction to the theory of hyperbolic buildings and their lattices, see the survey \cite{thomas} by A.~Thomas.

Here we study the 23 groups further, motivated by Gromov's famous surface subgroup question: Does every one-ended hyperbolic group contain a subgroup which is isomorphic to the fundamental group of a closed surface of genus at least 2? Recall, that a group $G$ is a {\em surface group} if $G=\pi_1(\mathcal{F})$, where $\mathcal{F}$ is a closed surface. If in addition $\mathcal{F}$ is finite, then $\pi_1(\mathcal{F})$ has one of the following forms (see \cite{stillwell}) 
\begin{enumerate}
 \item[(i)] $G=\langle a_1,b_1,\ldots, a_n,b_n\; |\: a_1b_1a_1^{-1}b_1^{-1}\cdots a_nb_na_n^{-1}b_n^{-1}\rangle$, when $\mathcal{F}$ is orientable and of genus $n$,
\item[(ii)] $G=\langle a_1,\ldots, a_n\; | \: a_1^2a_2^2\cdots a_n^2\rangle$, when $\mathcal{F}$ is non-orientable and of genus $n$.
\end{enumerate}
 
Even though Gromov's question remains open, it is known that there are many classes of hyperbolic groups, for which the answer is positive. For example, see \cite{crisp-sageev-sapir} for surface subgroups of right-angled Artin groups or \cite{calegari}, where Calegari and Walker show, that a random group contains many quasiconvex surface subgroups. Existence of surface subgroups in right-angled hyperbolic buildings
was shown in \cite{futer-thomas}, and in hyperbolic buildings with 4-gonal
apartments in \cite{Vdovina2005}. Existence of surface subgroups in fundamental groups of higher-dimensional complexes is discussed, for example, in \cite{reid}.

In \cite{KV3} the authors studied surface subgroups of the 23 torsion-free groups acting simply transitively on the vertices of hyperbolic triangular buildings of the smallest non-trivial thickness by searching periodic apartments of genus 2. They showed, that only in the first two groups listed in \cite{KV2} such an apartment exists, and therefore there also exists a surface subgroup of genus 2. Here we take another approach. Instead of studying the group on the hyperbolic building, we study first the 1-vertex polyhedron glued from the 15 triangles, and then the 3-cover of it, a 3-vertex polyhedron. We show that these polyhedra admit surface subgroups only in the same two cases out of 23, that were obtained in \cite{KV3}.

\section{Notations}

Let us consider the triangular complex $X$ glued together from 15 triangles with words on their boundary (for example as in Table \ref{T1T21}) . We relax on the condition of orientability, and consider the vector space $C_2$ freely spanned by the set of 
triangles of the triangular complex $X$, the ground field being 
$\frac{\integers}{2\integers}$. Thus the vectors in $C_2$ are sums of 
certain triangles (with coefficient 1). Now we also have the vector  
space $C_1$ freely spanned (on $\frac{\integers}{2\integers}$) by the 
set of edges of $X$. The boundary map $\partial : C_2\to 
C_1$ is defined on the basis of triangles by the relation :
$\partial T=a+b+c$ where $a,b,c$ are the three edges in the boundary 
of the triangle $T$. Thus, a 2-cycle (mod. 2) is an element of the kernel of $\partial$.

The triangular complex  $X$ is of thickness 3. This 
implies that any (non zero) 2-cycle corresponds to the immersion of a 
surface. Since we chose the ground field to be
$\frac{\integers}{2\integers}$ instead of $\integers$, the immersed surface may 
be non-orientable.

In order to really get  a (quasi-convex) surface subgroup we 
need also a geometric condition. One sufficient condition is that all local cycles are of 
length 8. We might allow a bit more, then the immersed surface would 
not necessarily be isometrically embedded anymore, only a quasi-isometric embedding.

We will divide the search of surface subgroups in the 23 groups into two cases: first we consider cycles in the 1-vertex polyhedron defined by the triangles, then in the 3-vertex cover of this polyhedron.

\section{Surfaces immersed to the 1-vertex polygon}

We calculated the 2-cycles (mod 2) in the one vertex case for all the 23 torsion free groups (that is, 23 sets of 15 triangles that each form the 1-vertex polygon, listed in the appendix II of \cite{KV2}) by going through all $2^{15}$ possible chains. It turned out that some of them are indeed cycles (mod 2). The cycles are listed in Table \ref{2cycles1vertex}. 
The notation "Presentation nb 1: 6-cycle 101101000000011"
means, that in the first polygonal presentation obtained in \cite{KV2}, there is one 2-cycle, it is of length 6, and it is formed by taking the 1st, 3rd, 4th, 6th, 14th, and 15th triangle in $T_1$ (see Table \ref{T1T21} for the triangles).

\begin{table}[h!]
\caption{2-cycles in the 1-vertex polyhedra}\label{2cycles1vertex}
\begin{center}
\begin{tabular}{| l | l | l | }
\hline
 Presentation nb           1
&6-cycle& 101101000000011\\

 Presentation nb           2
&6-cycle& 000011100010110\\
&6-cycle& 101101000000011\\
&8-cycle& 101110100010101\\

 Presentation nb           3
&8-cycle& 000011011110101\\

 Presentation nb           4
&8-cycle &000011011110101\\

 Presentation nb           5
&6-cycle& 000011001100011\\
&8-cycle& 101101010001110\\
&10-cycle& 101110011101101\\

 Presentation nb           6
&8-cycle& 001101110101100\\
&6-cycle& 100000001111001\\
&10-cycle& 101101111010101\\

 Presentation nb           7
&6-cycle& 000010110011100\\
&8-cycle& 001101110101100\\
&6-cycle& 001111000110000\\
&6-cycle& 100000001111001\\
&8-cycle& 100010111100101\\
&10-cycle& 101101111010101\\
&8-cycle& 101111001001001\\

 Presentation nb           8
&6-cycle& 101100010000011\\

 Presentation nb           9
&8-cycle& 000001111111010\\

 Presentation nb          10
&8-cycle& 100000110110111\\

 Presentation nb          11 
&8-cycle& 100001100111101\\

 Presentation nb          12 
&8-cycle& 100011000011111\\

 Presentation nb          13 
&6-cycle& 100011000111000\\

 Presentation nb          14 
&8-cycle& 100011011100110\\

 Presentation nb          15 
&10-cycle& 101110011010111\\

 Presentation nb          16 
&6-cycle& 100000110000111\\

 Presentation nb          18 
&6-cycle& 000000001111011\\

 Presentation nb          20 
&6-cycle& 001101001100100\\
\hline
\end{tabular}
\end{center}
\end{table}

To see whether the 2-cycles come from a surface subgroup, we need the geometric condition. To any cycle mod 2 there is a naturally associated compact triangulated surface together with an immersion inside the triangular complex. Each triangle is equipped with the regular hyperbolic metric with angle $\pi/4$ at vertices. Therefore if the initial cycle was locally an union of 8-cycles, it forces the immersion to be a local isometry, and it follows that the immersion is $\pi_1$-injective, because we are in the context of negative curvature.

Thus we know that if the 2-cycle intersects with the link along 1-cycles of length 8, it corresponds to a surface. This may be too strong condition, but at least it gives some examples.

The first candidate for a surface subgroup is in the presentation $T_9$: the 2-cycle of the triangles $(x_8,x_7,x_3)$, $(x_{10}, x_{13}, x_3)$, $(x_8,x_5,x_4)$, $(x_{14}, x_{14}, x_4)$, $(x_{10}, x_{12}, x_5)$, $(x_7,x_{12},x_6)$, ($x_{15}, x_9,x_6)$, $(x_{15},x_{13},x_9)$ has a vertex-link consisting of 3 circles of length 8. Thus the surface has 8 regular hyperbolic triangles with angles $\pi/4$, and around the vertex the triangles are situated in such a way that the geometry condition is fulfilled. A surface subgroup consisting of the same triangles was found in $T_9$ with a completely another method in also \cite{KV3}. 


Unfortunately, in all the other 22 cases the 2-cycles have vertex-links other than a set of cycles of length 8, and thus we cannot be sure of the existence of any other surface subgroups.

\section{Surfaces immersed to the 3-vertex cover}

Let us consider the 3-covers $Y_i$ of the polyhedra $X_i$ with one vertex. We construct the 3-covers as in \cite{CKV}: 
We substitute
each triple of the form  $(x_i,x_j,x_k)$ by three triplets
$(x_i^1,x_j^2,x_k^3)$, $(x_i^2,x_j^3,x_k^1)$ and
$(x_i^3,x_j^1,x_k^2)$, and each $(x_i,x_j,x_j)$ similarly by three
triplets $(x_i^1,x_j^2,x_j^3)$, $(x_i^2,x_j^3,x_j^1)$ and
$(x_i^3,x_j^1,x_j^2)$. We then have 45 triangles, which represent the
generators of the index 3 subgroup.
 Since the 1-vertex polygon $X_i$ consists of 15 triangles, the 3-cover $Y_i$ is glued together from 45 triangle faces. To find surface subgroups we first search for all 2-cycles (mod 2) in the 3-vertex polyhedron $Y_i$. 

In the 3-vertex space $Y$, we first pick some vertex $v$. Then all triangles of $Y$
  are adjacent to $v$. So $Y$ is obtained by identifying some boundary
edges of a complex $B$ which is a 1-ball around $v$: in other words $B$ is the star of a vertex $v_0$ in the universal
cover, so $B$ is the cone with base the generalized 4-gon $L$ and
apex $v$. We note that since each edge in $Y$ is contained in three triangles, the
identifications on the boundary $L$ of $B$ is very simple: it must
identify edges three by three. 


Now having a 2-cycle $c$ in $Y$ is the same as having a 2-chain $b$ in
$B$, whose boundary is contained in $L$ (so $b$ is locally a 2-cycle
at $v$), such that moreover the projection of $b$ is $c$ - so the free
faces of the 2-chain $b$ cancel pairwise. This is easy to express  only in terms of the graph $L$.
The 2-chain $b$ corresponds to a unique 1-chain $w$ in $L$. The
condition that $c$ is a 2-cycle exactly means that 
\begin{itemize}
\item $w$ is a 1-cycle: it's just a union of circles of the graph $L$,
\item the identification of the edges on $L$ identifies the edges of
$w$ two by two: when $e$ is an edge of $L$ appearing in $w$, it is
identified in $L$ with two other edges $e',e''$, and exactly one of
$e',e''$ belongs to $w$, the other one is exterior to $w$.
\end{itemize}
In order to find all possible such 1-cycles we
first write $L$ as the quotient of a rooted tree $T$ which is a ball
of radius 4 in the trivalent tree. So $L$ is obtained by an
identification of points at $\partial T$. All 1-cycles in the link are generated by 16 generators (cycles of length 8 in the link), easily visible from the rooted tree: each generator goes from one branch to another through the root. 

Then we calculate which of the $2^{16}$ combinations of these 16 generators give a 1-cycle $w$ which corresponds to a 2-cycle $b$ in the ball $B$, and thus to a 2-cycle $c$ in $Y$. 

With the help of a computerized algorithm, we find several such cycles for each of the 23 presentations.

Since we have a 2-cycle (mod 2) in $Y_i$, the cover of degree 3 of $X_i$, it means we have an immersion of a (not necessarily orientable) surface $S_i$ into $Y_i$. Composing with the covering map $p_i:Y_i \rightarrow X_i$ we get an immersion $i:S_i \rightarrow X_i$. Now we have to insure that we have a local isometry of non-positively curved complexes, i.e. the map is a $\pi_1$-injection, to guarantee that we have found surface subgroup acting on the building.

In this context a local isometry of a triangular surface $S_i$ to the 
triangular complex $X_i$ is a locally injective simplicial map $f_i:S_i\rightarrow  X_i$, 
such that the image of the circle link($S_i,v$) is a circle of link($X_i,f_i(v)$) 
with no short cut of length $<\pi$, when each edge of the link receives the 
length $\pi/4$. So the condition means that if there is a path of 
less than $4$ edges between two vertices of $f_i(link(S_i,v))$ then this path is 
inside the circle $f_i(link(S_i,v))$.

%
%
%


When we go through all the 23 torsion free groups that act simply transitively on the vertices of hyperbolic triangular buildings of the smallest non-trivial thickness, we find 2-cycles in most of them. In six cases the link at least at one of the three vertices is a collection of cycles of length 8. These are presented in the Appendix. In three of the cases the links at all vertices are cycles of length 8. One of them is the same candidate for a surface subgroups in the group $G_9$ that we found already in the 1-vertex polyhedron.
The other two cases, were surface subgroups are four, are the same two groups $G_1$ and $G_2$ for which the authors in \cite{KV3} find periodic planes of genus 2 and therefore also surface subgroups of genus 2.

\subsubsection*{Group $G_1$}

We find a 2-cycle in the 3-vertex cover. It consists of 8 triangles: $(x_1^3, x_{10}^2, x_1^1)$, $(x_{15}^3, x_2^2, x_1^1)$, $(x_{15}^2, x_2^1, x_1^3)$, $(x_{11}^3, x_9^2, x_2^1)$, $(x_{14}^1, x_3^3, x_2^2)$, $(x_{15}^2, x_{13}^1, x_3^3)$, $(x_{14}^1, x_{15}^3, x_9^2)$ and $(x_{13}^1, x_{11}^3, x_{10}^2)$. On the 3-vertex polygon, the link at each vertex is a circle of length 8. On the 1-vertex polygon the image of this 2-cycle under the covering map has therefore a link which is a collection of 3 cycles of length 8. Since the complete vertex-link is a generalized 4-gon and has girth 8, there cannot be any shortcuts  between the vertices of the cycles. Thus we have a local isometry, and therefore there exists a surface subgroup in the group $G_1$. 

 This surface subgroup is not the only one, with the method we find two other surface subgroups in $G_1$ as well. These too come from 8 triangles each, and the links at each vertex are cycles of length 8.


\subsubsection*{Group $G_2$}

Also in $G_2$ we find a surface subgroup. It comes from the same triangles $(x_1^3, x_{10}^2, x_1^1)$, $(x_{15}^3, x_2^2, x_1^1)$, $(x_{15}^2, x_2^1, x_1^3)$, $(x_{11}^3, x_9^2, x_2^1)$, $(x_{14}^1, x_3^3, x_2^2)$, $(x_{15}^2, x_{13}^1, x_3^3)$, $(x_{14}^1, x_{15}^3, x_9^2)$ and $(x_{13}^1, x_{11}^3, x_{10}^2)$ in the 3-cover as in the previous case. This is sure to be a surface, since the links at the vertices are cycles of length 8. Again, this is not unique. There are two other sets of 8 triangles in $T_2$ that give a surface subgroup.

\subsubsection*{Group $G_9$}

For $G_9$ we already found a surface from the 1-vertex polygon. The corresponding surface is naturally found from the 3-cover, too. On the 3-cover the surface consists of 24 triangles, since every triangle from the 1-vertex polygon is presented three times here. 

\subsubsection*{Group $G_{14}$}
24 triangles in the 3-cover give at one of the vertices a link that consists of 3 cycles of length 8, but at another vertex a cycle of length 24. Such a long cycle always has a shortcut somewhere: take two points from the link that have in this cycle four other points between them. Since the diameter of the generalized 4-gon is four, there is always a path of length four between these points. Thus the local isometry is violated and this 2-cycle does not correspond to a surface subgroup.

\subsubsection*{Group $G_{18}$}
16 triangles in the 3-cover form a cycle, but even though at two of the vertices the links consists of 2 cycles of length 8, at the third vertex it is just one cycle of length 16. Thus, similarly as above, this 2-cycle does not correspond to a surface subgroup.

\subsubsection*{Group $G_{20}$}
24 triangles in the 3-cover form a cycle, but at one of the vertices the link is a cycle of length 24, at another it consists of cycles of lengths 16 and 8, and only in one vertex the link is a combination of cycles of length 8. In the long cycles there are definitely some shortcuts, so this 2-cycle does not correspond to a surface subgroup.

\section*{Acknowledgements}
The author is extremely grateful for all the discussions and support from Fr\'ed\'eric Haglund and Alina Vdovina.


\newpage 

\section*{Appendix: Cycles in 3-vertex polyhedra}

Here the labels $x_i^j$, $i=1, \ldots, 15$, $j=1,2,3$, of the sides of the triangles are given as integers from 1 to 45, $x_i^j\hat{=}3(i-1)+j$.

\begin{verbatim}
 Number           1

 Triangles in the 2-cycle:
  3 29  1
 45  5  1
 44  4  3
 33 26  4
 40  9  5
 44 37  9
 40 45 26
 37 33 29
 The cycles are generated by:
  9 40 45  1  3  4 33 37
  Cycles in the link:
  1  3  4 33 37  9 40 45  
  
 Triangles in the 2-cycle:
  2 28  3
 43  6  2
 44  4  3
 32 25  6
 42  8  4
 43 39  8
 42 44 25
 39 32 28
 The cycles are generated by:
  6  7 45  1  3 28 39 43
  6  7 45  1  3  4 42 25
  Cycles in the link:
  3  4 42 25  6 43 39 28  

 Triangles in the 2-cycle:
  1 30  2
 45  5  1
 43  6  2
 31 27  5
 41  7  6
 45 38  7
 41 43 27
 38 31 30
 The cycles are generated by:
 27 31 30  1  3 28 39 43
  6  7 45  1  3 28 39 43
  Cycles in the link:
  1 30 31 27 43  6  7 45  


 Number           2

 Triangles in the 2-cycle:
  3 29  1
 45  5  1
 44  4  3
 33 26  4
 40  9  5
 44 37  9
 40 45 26
 37 33 29
 The cycles are generated by:
  9 40 45  1  3  4 33 37
  Cycles in the link:
  1  3  4 33 37  9 40 45  

 Triangles in the 2-cycle:
  2 28  3
 43  6  2
 44  4  3
 32 25  6
 42  8  4
 43 39  8
 42 44 25
 39 32 28
 The cycles are generated by:
  6  7 45  1  3 28 39 43
  6  7 45  1  3  4 42 25
  Cycles in the link:
  3  4 42 25  6 43 39 28  

 Triangles in the 2-cycle:
  1 30  2
 45  5  1
 43  6  2
 31 27  5
 41  7  6
 45 38  7
 41 43 27
 38 31 30
 The cycles are generated by:
 27 31 30  1  3 28 39 43
  6  7 45  1  3 28 39 43
  Cycles in the link:
  1 30 31 27 43  6  7 45  


 Number           9

 Triangles in the 2-cycle:
 24 20  7
 22 21  8
 23 19  9
 30 38  7
 28 39  8
 29 37  9
 24 14 10
 22 15 11
 23 13 12
 42 41 10
 40 42 11
 41 40 12
 30 35 13
 28 36 14
 29 34 15
 21 35 16
 19 36 17
 20 34 18
 45 26 16
 43 27 17
 44 25 18
 45 38 25
 43 39 26
 44 37 27
 The cycles are generated by:
 15  7 30  1  3 28 39 34
 15  7 30  1  3  4 33 22
 24  7 30  1  3  4 42 10
 27 13 30  1  3 28 39 43
 27 13 30  1  3  4 33 37
 12 13 30  1  3  4 42 40
 21 16 45  1  3  4 33 22
  9 25 45  1  3 28 36 19
  9 25 45  1  3  4 33 37
 18 25 45  1  3 28 39 34
  Cycles in the link:
  7 24 10 42 40 12 13 30  
  9 19 36 28 39 43 27 37  
 15 22 21 16 45 25 18 34  


 Number		14

 Triangles in the 2-cycle:
  2 43  3
 28  6  2
 29  4  3
 33  8  4
 32  7  6
 20 10  9
 45 35  7
 43 36  8
 44 34  9
 24 17 10
 22 18 11
 34 27 11
 22 39 14
 40 24 14
 21 26 16
 19 27 17
 36 32 16
 35 31 18
 30 38 19
 28 39 20
 29 37 21
 40 30 26
 45 38 31
 44 37 33
 The cycles are generated by:
  6  7 45  1  3 43 39 28
 21  7 45  1  3 43 36 16
 21  7 45  1  3  4 33 37
  9 31 45  1  3 43 36 10
  9 31 45  1  3  4 33 34
 18 31 45  1  3 43 39 22
 27 19 30  1  3  4 33 34
 24 40 30  1  3 43 36 10
  Cycles in the link:
  3  4 33 37 21 16 36 43  
  6  7 45 31 18 22 39 28  
  9 10 24 40 30 19 27 34  
 



 Number          18

 Triangles in the 2-cycle:
 33 11  4
 32 10  6
 42  8  4
 41  7  6
 27 14  7
 25 15  8
 42 23 10
 40 24 11
 17 34 15
 43 39 14
 19 27 17
 33 29 19
 40 36 23
 41 34 24
 39 32 25
 43 36 29
 The cycles are generated by:
  6 28 45  1  3 43 39  7
  6 28 45  1  3  4 42 10
 27 31 30  1  3 43 39  7
 27 31 30  1  3  4 33 19
 15 34 30  1  3 43 39 25
 24 34 30  1  3 43 36 40
  Cycles in the link:
  4 33 19 27  7  6 10 42  
 15 25 39 43 36 40 24 34  

 Triangles in the 2-cycle:
 31 12  5
 32 10  6
 40  9  5
 41  7  6
 27 14  7
 26 13  9
 42 23 10
 41 22 12
 16 36 14
 45 38 13
 21 26 16
 32 28 21
 42 35 22
 40 36 23
 38 31 27
 45 35 28
 The cycles are generated by:
  9 13 45  1  3 43 36 40
  6 28 45  1  3 43 39  7
  6 28 45  1  3  4 42 10
 21 28 45  1  3 43 36 16
 12 31 30  1  3  4 42 22
 27 31 30  1  3 43 39  7
  Cycles in the link:
  6  7 27 31 12 22 42 10  
  9 13 45 28 21 16 36 40  


 Number          20

 Triangles in the 2-cycle:
  3 29  1
  2 28  3
 45 17  1
 43 18  2
  9 20  4
  7 21  5
 36 23  4
 34 24  5
 13 12  8
 14 10  9
 33 41  7
 31 42  8
 45 38 10
 44 37 12
 27 44 13
 26 43 15
 28 36 14
 29 34 15
 40 33 17
 41 31 18
 24 38 19
 23 37 21
 42 26 19
 40 27 20
 The cycles are generated by:
 12 13 30  1  3 28 36 37
 27 13 30  1  3 16 33 40
 24 34 30  1  3 16 42 19
 15 34 30  1  3 28 39 43
  9 10 45  1  3 28 36  4
 18 10 45  1  3 28 39 43
 18 10 45  1  3 16 42 31
 21 25 45  1  3 28 36 37
 21 25 45  1  3 16 33  7
  Cycles in the link:
  1  3 28 36  4  9 10 45  
  7 33 40 27 13 12 37 21  
 15 34 24 19 42 31 18 43  

\end{verbatim}

\end{document}